\documentclass[12pt,reqno]{amsart}

\usepackage{amssymb,amsmath,eucal}%,epsf,graphicx,epic}

\renewcommand{\d}{\delta}

\newcommand{\f}{\varphi}

\newcommand{\s}{\sigma}

\newcommand{\G}{\Gamma}

\newcommand{\La}{\Lambda}

\font\cmssl=cmss10 at 12 pt

\newcommand{\bR}{\mathbb{R}}
\newcommand{\bZ}{\mathbb{Z}}

\newcommand{\cC}{\mathcal{C}}

\newcommand{\cI}{\mathcal{I}}

\newcommand{\cP}{\mathcal{P}}

\DeclareMathOperator\vol{vol}

\DeclareMathOperator{\Gr}{Gr}

\newtheorem{Th}{Theorem}
\newtheorem{Prop}{Proposition}
\newtheorem{Cor}{Corollary}
\newtheorem{Lem}{Lemma}
\newtheorem{Def}{Definition}
\newtheorem{Ex}{Example}

\newcommand{\bpr}{\begin{Pres}\ \ }
\newcommand{\epr}{\end{Pres}\ \ }
\newcommand{\bex}{\begin{Ex}\ \ }
\newcommand{\eex}{\end{Ex}\ \ }
\newcommand{\bt}{\begin{Th}\ \ }
\newcommand{\et}{\end{Th}\noindent}
\newcommand{\bp}{\begin{Prop}\ \ }
\newcommand{\ep}{\end{Prop}\noindent}
\newcommand{\bc}{\begin{Cor}\ \ }
\newcommand{\ec}{\end{Cor}}
\newcommand{\bl}{\begin{Lem}\ \ }
\newcommand{\el}{\end{Lem}}
\newcommand{\bd}{\begin{Def}\ \ }
\newcommand{\ed}{\end{Def}\noindent}
\newcommand{\pf}{\noindent{\it Proof:\ \ }}
\newcommand{\be}{\begin{equation}}
\newcommand{\ee}{\end{equation}}
\newcommand\la[1]{\label{#1}}
\newcommand\re[1]{(\ref{#1})}
\def\<#1,#2>{\langle\,#1,\,#2\,\rangle}
\newcommand{\arr}{\begin{array}{rlll}}
\newcommand{\ea}{\end{array}}
\newcommand{\bea}{\begin{eqnarray}}
\newcommand{\eea}{\end{eqnarray}}
\newcommand{\bean}{\begin{eqnarray*}}
\newcommand{\eean}{\end{eqnarray*}}

\newcommand{\1}{   1 \hspace*{-0.87ex}
                   1 }%\hspace*{0.9ex}}

\begin{document}
\rightline{math.DG/0604558}
\vskip 1 cm
\title[Special Graphs]{Special Graphs}
\author[C. Devchand, J. Nuyts]{Chandrashekar Devchand, Jean Nuyts} 
\address{Mathematisches Institut\\
Universit\"at Bonn\\
Beringstr. 1\\
D-53115 Bonn}
\email{devchand@math.uni-bonn.de}
\email{gw@sas.upenn.edu}

\author[G. Weingart]{Gregor Weingart} 
\address{%Physique Th\'eorique et Math\'ematique\\
Universit\'e de Mons-Hainaut\\ 20 Place du Parc\\ B-7000 Mons\\ Belgium}
\email{Jean.Nuyts@umh.ac.be}
\dedicatory{Dedicated to Dmitri V. Alekseevsky on the occasion of his sixty-fifth birthday}
 \thanks{
This work was supported by the {\it Schwerpunktprogramm Stringtheorie} of the
Deutsche For\-schungs\-ge\-mein\-schaft, MPI f\"ur Mathematik, Bonn
and the Belgian {\it Fonds National de la Recherche Scientifique}}

\begin{abstract}
A {\em special} $p$-form is a $p$-form which, in some orthonormal basis $\{e_\mu \}$, has components 
$\f_{\mu_1\ldots\mu_p} = \f(e_{\mu_1},\ldots, e_{\mu_p})$ taking values in $\{-1,0,1\}$. 
We discuss graphs which characterise such forms.  
\end{abstract}
\maketitle

\section{Calibrations, special forms and graphs}

A constant $p$-form $\f$ in a $d$-dimensional
Euclidean space is a {\it calibration} if for any $p$-dimensional subspace
spanned by a set of orthonormalised vectors $e_1,\dots, e_p$, the following condition
holds:
\begin{equation}
(\f(e_1,\dots,e_p))^2 \leq 1\ ,
\end{equation}
with equality holding for at least one subspace.
Let $U$ be an oriented $p$--dimensional subspace of $\bR^d$
with oriented metric volume $\vol_U$. The set of all such subspaces is the oriented
Grassmannian $\Gr_p\,\bR^d$. 
A calibration $\varphi\,\in\,\La^p\,\bR^d$  is thus a $p$--form
with the property that the function 
$\overline\varphi:   \Gr_p\,\bR^d\;\longrightarrow\;\bR$ associated to $\varphi$
and defined by
$U \mapsto \overline\varphi(U):=\langle\varphi,\vol_U\rangle$ 
takes values in $[-1,1]\subset\bR$ with at least one of the two extremal values $\pm 1$
being achieved. 
The $p$-planes $U$ for which $\overline\varphi(U)=\pm 1$ are said to be calibrated by $\varphi$.

Almost all examples of calibrations known are invariant
under a group $G\subset {\rm O}(\bR^d)$ large enough 
so that it is relatively simple to check the calibration condition directly.
Interestingly most of these examples, in particular the calibrations characterising
special holonomy manifolds, for instance the $G_2$--invariant Cayley 3-form in seven
dimensions, defined by the structure constants of the imaginary octonions, and 
the Spin(7)-invariant 4-forms in eight dimensions are special forms:
 
\bd\la{special}
 A {\cmssl special p--form} $\f$ is a $p$--form $\f\,\in\,\La^p\bR^d$ on
 $d$--dimensional Euclidian space $\bR^d$ in the orbit under the orthogonal 
 group $\mathrm{O}(d,\bR)$ of
 \begin{equation}
  \f\;\;=\;\;\sum_{1\leq \mu_1<\ldots<\mu_p\leq d}
  \f_{\mu_1\ldots\mu_p}\,e_{\mu_1}\wedge e_{\mu_2}\wedge\ldots\wedge e_{\mu_p}
 \la{special1}
 \end{equation}
 with $\f_{\mu_1\ldots\mu_p}\,\in\,\{-1,0,1\}$ and  
 $(e_1,\dots, e_d)$ an orthonormal basis. 
\ed
In other words, a $p$-form $\f$ is special if there exist $d$ orthonormal basis vectors 
$ e_\mu, \,\mu=1,\ldots,d$, such that for any subset of $p$ basis vectors
$e_{\mu_1},\ldots, e_{\mu_p}$ we have
\begin{equation}
 \f_{\mu_1\ldots\mu_p} \;\;:=\;\;
 \f(e_{\mu_1},\ldots, e_{\mu_p}) \;\;\in\;\; \{-1,0,1\}\ .
 \la{special2}
\end{equation}
Given a basis $\{e_{\mu}\}$, there are clearly only a finite number (obviously less than 
$3^{\frac{d!}{p!(d-p)!}}$) of orbits of special $p$--forms under
$\mathrm{O}(d,\bR)$ parametrised by the components $\f_{\mu_1\ldots\mu_p}\,
\in \,\{-1,0,1\}$.
Apparently different special $p$--forms may nevertheless 
be in the same orbit under $\mathrm{O}(d,\bR)$, because the
subgroup $\mathrm{O}(d,\bZ)\,\subset\,\mathrm{O}(d,\bR)$ of orthogonal matrices
with integer coefficients maps the special form $\f$ in equation \re{special1}
again into a special form with possibly different components. The group $\mathrm{O}(d,\bZ)$ 
is isomorphic to the semidirect product of the permutation group acting naturally on 
$d$ copies of $\bZ_2$. 
The action of 
$(\s,\eta_1,\dots,\eta_d)\,\in\,S_d\ltimes\bZ_2^d\,\cong\, \mathrm{O}(d,\bZ)$
on the antisymmetric tensor indices of $\f$ is given by 
$\f_{i_1\ \dots\ i_p} \mapsto \eta_{i_1}\dots\eta_{i_p}\,\f_{\s(i_1)\ \dots\ \s(i_p)}$, 
where $\s\in S_d$ and $\eta_i^2=1\,,\,i=1,\dots,d$.

Let us now give an alternative description of special forms. An oriented 
$p$-subset of $\{1,2,\dots,d\}$ is given by the $p$ elements $s=\{\mu_1,\dots,\mu_p\}$ 
such that $1\le \mu_1 < \mu_2 < \dots < \mu_p\le d $. The space of all such $p$-subsets 
is the {\em vertex space} $\cP^p(\{1,\dots,d\})$. These oriented subsets of 
$\{1,2,\dots,d\}$ are in bijective correspondence to oriented coordinate subspaces 
$\bR^p\,\subset\,\bR^d$ via 
$s=\{\mu_1,\ldots,\mu_p\}\,\longmapsto\,e_{\mu_1}\wedge\ldots\wedge e_{\mu_p}$. 
A special $p$--form can be thought of as a function from $\cP^p(\{1,\dots,d\})$
to $\f_{\mu_1\ldots\mu_p}\,\in\,\{-1,0,1\}$.
Consequently a special $p$--form is specified completely by the
two sets $\cI^+$ and $\cI^-$ of oriented subsets $\{\mu_1,\dots,\mu_p\}\,\subset\,
\{1,\dots,d\}$ which have respectively $\f_{\mu_1\ldots\mu_p}\,=\,+1$ and $-1$. 
We denote by $\mu^{(a)},a=1,\dots,|\f|$, the elements in $\cI\,:=\,\cI^+\cup\cI^-$,
the support of $\f$ and
\begin{equation}
\f = \sum_{a=1}^{|\f|} \f_{\mu_1^{(a)}\dots \mu_p^{(a)}} \ 
                       e_{\mu_1^{(a)}}\wedge \dots\wedge e_{\mu_p^{(a)}}\ ,
\la{mua}
\end{equation}
where $|\f|\,:=\,|\cI|$ is the weight of $\f$. 
For every permutation $\sigma\in S_d$ of the basis vectors $e_1,\dots,e_d$, there exists 
a corresponding permutation of the oriented $p$-subsets.

Interestingly, we can define a metric on $\cP^p(\{1,\dots,d\})$ by setting 
the distance between two oriented $p$-subsets, $s$ and $\tilde s$, to be
$
\d(s,\tilde s) = p - \#(s \cap \tilde s)
$,
where $\#(s \cap \tilde s)$ is the number of elements in the intersection 
of the sets $s$ and $\tilde s$.
We can visualise
the restriction of this metric to the set $\cI$ by drawing a graph with
labeled edges, the vertices corresponding to the elements of $\cI$ and the
edges running between vertices labeled by a {\em distance} strictly less than 
$p$. Unfortunately the graph of a special $p$--form
$\f$ does not specify the components
$\f_{\mu_1\dots\mu_p}\,\in\,\{-1,0,1\}$ completely up to
the action of $\mathrm{O}(d,\bZ)$; we still need to specify some relative
sign. Nevertheless the graph gives a very condensed way of encoding the
characteristics of a special $p$--form. 

Consider a graph $\G$ composed of a set of vertices $V=\{v_i\ ;\ i=1,\dots r \}$ 
connected by the maximum possible number of edges, $r(r-1)/2$, each
labeled by a positive number $d(v_i,v_j)$, the {\em distance} between the vertices
$v_i$ and $v_j$ at its ends.
A graph is {\em admissible} if any triangle with edges labeled by distances
$d_i,d_j,d_k$ satisfies the triangle inequalities
\be
1\le d_i \le d_j+d_k  \quad {\rm{and\ cyclic\ permutations}}.
\end{equation}

\bd
A {\em realisation} of a graph $\G$ is a map
\bea
\rho: \G &\rightarrow& \cP^p(\{1,\dots,d\}) 
         \nonumber\\
    v &\mapsto& s_v
\eea
which assigns to any vertex $v$ an oriented $p$-subset $s_v$ 
such that the distance between any two
vertices $d(v,w)$ is equal to the distance between the corresponding
oriented $p$-subsets $\d(s_v,s_w)$ and that 
\bea
\# \left( \bigcup_{v\in \G} s_v \right) = d \quad ,\quad
\# \left( \bigcap_{v\in \G} s_v \right) =0\ .
\eea
\ed
Two realisations are equivalent if there exist a permutation $\s\in S_d$
of the oriented $p$-subsets which maps one onto the other.

Consider the power set $\cP(V)$, the set of all subsets $S\subset V$ of 
vertices of the graph $\G$,
\begin{equation}
\cP(V) = \{\{\emptyset\},\{v\},\{w\},\dots,\{v,w\},\{x,y\},\dots,\{v,w,x\},\dots,\{V\}\}\ .
\end{equation}
Clearly, $\#\{\cP(V)\}=2^r$. 
For every realisation, a graph function $f$ associates 
a nonnegative integer to every $S\in \cP(V)$ as follows
\bea
F(S) &:=& \left( \left( \bigcap_{v\in S} s_v\right) \bigcap  
                  \left( \bigcup_{v\notin S} s_v \right)^{\!\!\!\! {\rm{C}}}  \right)\ge 0
           \nonumber\\
f(S) &:=& \# (F(S))\  ,
\eea
where C denotes the complementary subset.
Clearly, this function measures the number of indices which occur in every
$s_v$ for $v\in S$ but do not occur in any $s_v$ for $v\notin S$. Trivially,
we have
\bea
f(\emptyset)&=&0 
       \nonumber \\  
f(V)&=&0    
       \nonumber \\
\sum_{S\in \cP(V)}f(S)&=&d\ ,
\label{bc}
\eea 
because ${\displaystyle{\bigcap_{v\in S} s_v}}$ is empty for the two first cases 
and $\displaystyle{\sum_{S\in \cP(V)} f(S)}$ contains all the indices $\{1\dots d\}$ and 
each index contributes to $f$ for one and only one subset $S$.

\bt
To every graph function $f$, with non negative integer values, which satisfies 
\be
d(v,\tilde v)=p \quad - \sum_{\{    S\in \cP(V)\mid     
                            v,      
                            \tilde v\in S  \} }   f(S)\ ,
\label{basiceq}
\end{equation}
there corresponds a class of equivalent realisations of the graph.
\et
In particular, if $v=\tilde v$ 
\be
d(v, v)=p \quad - \sum_{\{S\in \cP(V) \vert
                  v\in S  \}  }   f(S)=0 \ .
\end{equation}
As a consequence, all realisations of a graph are simply obtained by finding all the 
non negative solutions to \re{basiceq}.
Every realisation yields a simple and direct construction of a special form, up to choices of signs.    
Examples can easily be generated \cite{dnw2}.

\section{Democratic Graphs}

Consider a graph $\G$ with vertices $v_i,i=1,\dots,r$, and the set of nonzero distances
%$\{d_a\ ;\ a=1,\dots,r(r{-}1)/2\}=
$\{d(v_i,v_j)\ ;\  i<j\}$.
The $r{\times}r$ {\em distance matrix} 
\be 
M_{ij}^{[r]}=d(v_i,v_j) 
\end{equation}
is clearly symmetric, with diagonal elements equal to $0$.

\bd 
A {\em symmetry} $\s$ of a graph $\G$ is a permutation of the vertices 
$v_i\mapsto \tilde v_i$ which leaves the distance matrix invariant, i.e.
\be
d(\tilde v_i,\tilde v_j)=d(v_i,v_j)\ . 
\end{equation}
\ed
If there exists a realisation of a graph with symmetry $\s$, such that $f(\s S)=f(S)$ for any $S\in\cP(V)$,
then the realisation has a permutation of the 
indices which induces $\s$ on the monomials.

\bd
A graph is {\em democratic} if, for every pair of vertices $v_i,v_j$, 
there  exists some symmetry $\sigma$ which maps 
$v_i\mapsto \tilde v_i=v_j$. 
\ed
Let $\{d_a\}$ denote the set of unequal distances.
\bp
A necessary condition for a graph with $r$ vertices to be democratic is that 
$n^{(i)}_a$, the number of vertices at distance $d_a$ from vertex $v_i$,
is independent of the choice of $v_i$. A graph satisfying this condition,
i.e. $n^{(i)}_a = n_a$, will be called {\em predemocratic}.
\ep
Note that $\sum_a n^{(i)}_a=r-1 $ and that a predemocratic graph is not necessarily democratic.

First, consider graphs with an even number of vertices, $r=2n$. 
Then, from every vertex there are $r{-}1$ edges labeled by $r{-}1$ distances.
To have democracy, every vertex should have the same set of distances to
its neighbouring vertices. Let $d_i:=d(v_1,v_{i+1}),\,i=1,\dots,r{-}1$, 
be the distances between $v_1$ and $v_{i+1}$. An example of a predemocratic 
graph with an even number of vertices $r=2n$ has distance matrix, 
up to relabeling, of the form
\bea
d(v_i,v_j)&=&(1-\delta_{ij})\ d_{\,j+i-2 \!\!\!\pmod{r-1}}
             \nonumber\\
d(v_i,v_{r})&=&(1-\delta_{ir})\ d_{\,2i-2 \!\!\!\pmod{r-1}}          
\eea  
with $d_0\equiv d_{r-1}$. 
For $r=4$, the corresponding graph is the unique predemocratic graph 
and it is also democratic. For higher $r$'s, these graphs are in general not democratic.

Now, consider graphs with an odd number of vertices, $r=2n{+}1$.

\bl\label{evenn}
If the number of vertices $r$ is odd, a predemocratic graph has
all  $n_a$'s even. 
\el
\pf
For a predemocratic graph with $r$ vertices, the total number of edges of
distance $d_a$ is clearly $n_a r/2$. \hfill $\square$

\noindent For $r$ odd, if we set $n_a=2$, for all $a$, 
there are $(r{-}1)/2$ unequal distances $d_a$.  
Distance matrices with $n_a=4,6,\dots$, with all distances $d_a$ different, can always be obtained 
from the distance matrices with $n_a=2$ by setting some $d_a$'s to be equal.

We now classify all distance matrices with  $n_a=2$ for all $a$ and all distances $d_a$ different.
Under these assumptions, the edges of length $d_a$ for given $a$ form a closed (possibly disconnected) 
curve $\cC_a$ containing every vertex once. Let us call the number of vertices 
in a connected piece of curve $\cC_a$ its pathlength, which is obviously
between 3 and $r$.

\bl
A necessary condition for a predemocratic graph with an odd number of 
vertices $r$ and $n_a=2$ to be democratic is that the curve $\cC_a$, for every $a$, 
has connected pieces of equal pathlength $3\le L_a \le r$. In other words, 
$L_a$ is a divisor of $r$ and $\cC_a$ consists of $r/L_a$ disconnected pieces.    
\el
\pf
The proof follows from the definition of a democratic graph:
the curve $\cC_a$ as seen from any vertex has the same form, 
independent of the choice of the vertex. 

\hfill$\square$

\noindent
To give an example of a democratic graph with $r=2n+1$ vertices
and $n_a=2$ for all $a$ we choose $n$ distinct positive integers
$d_1,\ldots,d_n$ and define the distance matrix $M^{[r]}$ by
\be
 M^{[r]}_{ii}=0
 \qquad,\qquad
 M^{[r]}_{ij}=d_{\mathrm{min}{\{|i-j|,2n+1-|i-j|\}}}
\label{pat1}\end{equation}
for $0\leq i,j\leq 2n$. Evidently the matrix $M^{[r]}$ has a cyclic
isometry group $\bZ_r$ shifting the vertices $v_i\longmapsto
v_{i+1(\mathrm{mod}\,r)}$. Assuming that the distances $d_1,
\ldots,d_n$ can be chosen in such a way that there exists a
graph function $f$ for the $M^{[r]}$ invariant under $\bZ_r$
we get a realisation for $M^{[r]}$, which is a democratic
graph with symmetry group containing $\bZ_r$.

\bt
 For an odd prime number $r=2n+1$ every democratic graph with $r$ vertices
 satisfying $n_a=2$ for all $a$ (and $n$ distinct distances $d_a$)
 has a distance matrix of the form $M^{[r]}_{ij}$ with a suitable
 choice of the positive integers $d_1,\ldots,d_n$.
\et

\pf
Essentially we only need to show that the symmetry group of a democratic
graph with a prime number $r=2n+1$ of vertices and $n_a=2$ for all $a$
must contain a cyclic subgroup of order $r$ acting transitively on the
vertices. Clearly all curves $\cC_a$ must be connected circles of length
$r$, because $r$ being prime has no proper divisors. We fix two vertices
$v_0$ and $v_1$ and the curve $\cC_a$ containing the edge between them.
By democracy there exists a symmetry $\s$ sending $v_0$ to $\s(v_0)=v_1$
and mapping the curve $\cC_a$ to itself, because it is a symmetry and all
$d_a$ are distinct. Thus $\s$ can only be the cyclic shift by one step
along the curve $\cC_a$, which clearly generates a cyclic group of
symmetries of order $r$ acting transitively on the vertices. Setting
$v_i:=\s^i(v_0)$ for $0\leq i\leq 2n$ we conclude that $M^{[r]}_{ij}:=
d(v_i,v_j)=d(v_0,v_{j-i})=d(v_{2n+1-j+i},v_0)$ for all $0\leq i\leq j\leq 2n$.
\hfill $\square$

Slightly more generally we can consider graphs with a group
of symmetries isomorphic to $\bZ_{r_1}\times\bZ_{r_2}$ acting 
transitively on the vertices. Of course if $r_1$ and $r_2$ are
relatively prime, then the group of symmetries considered is
isomorphic to $\bZ_{r_1r_2}$. Nevertheless we expect new features
compared to the classification above, because $r_1r_2$ is no longer
prime. With a group $\bZ_{r_1}\times\bZ_{r_2}$ acting 
transitively on the vertices it is convenient to label the
vertices by tuples $(i_1,i_2)\in\{1,\ldots,r_1\}\times\{1,\ldots,r_2\}$.
Straightening this out by replacing $(i_1,i_2)$ with $i:=i_2+r_2(i_1-1)$
we get a distance matrix of the form
\be
M^{[r_1][r_2]}=
\begin{pmatrix}
 M^{[r_2]} &Q_1       &Q_2       &Q_3       &\dots  &Q_2^t    &Q_1^t    \cr
 Q_1^t     &M^{[r_2]} &Q_1       &Q_2       &\ddots &Q_3^t    &Q_2^t    \cr
 Q_2^t     &Q_1^t     &M^{[r_2]} &Q_1       &\ddots &Q_4^t    &Q_3^t    \cr
 Q_3^t     &Q_2^t     &Q_1^t     &M^{[r_2]} &\ddots &Q_5^t    &Q_4^t    \cr
 \vdots    &\ddots    &\ddots    &\ddots    &\ddots &\ddots   &\ddots   \cr
 Q_2       &Q_3       &Q_4       &Q_5       &\ddots &M^{[r_2]}&Q_1      \cr
 Q_1       &Q_2       &Q_3       &Q_4       &\ddots &Q_1^t    &M^{[r_2]}\cr
\end{pmatrix}\quad ,
\label{M2}
\end{equation}
where $M^{[r_2]}$ is the $r_2\times r_2$ distance matrix defined in \re{pat1} 
which depends on $(r_2-1)/2$ arbitrary distances. The $r_2\times r_2$  
matrices $Q_i,\, i=1,\dots,(r_1{-}1)/2\, ,$ depend on $r_2$ arbitrary 
parameters. The first row of every $Q_i$ is arbitrary and the $r_2{-}1$ 
following rows are obtained by cyclically permuting the elements of 
the first row.  

Clearly a graph function $f$ for a matrix of the form \re{M2}, 
with the property that $f(\s S) = f(S)$ for all sets of vertices $S$ and all
$\s\in \bZ_{r_1}\times\bZ_{r_2}$, defines an equivalence class of democratic graphs.  
Conversely, for a democratic graph with $r=r_1r_2$ vertices, $r_1,r_2$ prime and $n_a=2$ for all $a$,
the distance matrix must be of the form  $M^{[r]}$ 
in  \re{pat1} or  $M^{[r_1][r_2]}$ in \re{M2}.

In full generality, for every factorisation $r=r_1r_2\cdots r_k$ with  
$r_1\ge r_2\ge  \dots\ge r_k > 1 $,  we can consider graphs having 
a group of symmetries isomorphic to $\bZ_{r_1}\times\dots\times\bZ_{r_k}$ 
acting transitively on the $r$ vertices. 
A convenient labeling of the
vertices is given by 
$v^{i_1,\dots, i_k} := \s_1^{i_1}\s_2^{i_2}\dots\s_k^{i_k} (v_0)$
where the tuples 
${i_1,\dots, i_k} \in \{1,\ldots,r_1\}\times\dots\times\{1,\ldots,r_k\}$ 
and the $\sigma_A , A=1,\dots,k$, are generators of  $\bZ_{r_{\!A}}$,
cyclic permutations of order $r_{\!A}$.
Then the distance matrices have elements
\be
 d(v^{i_1,\dots, i_k} , v^{j_1,\dots, j_k})
 = d(v_0 , v^{j_1-i_1,\dots, j_k-i_k})
 = d_{j_1-i_1,\dots, j_k-i_k}  \quad,\quad i_A < j_A\, ,\, A=1,\dots,k.
 \end{equation}
If we have a graph function $f$ for this matrix, with 
the property that $f(\s S) = f(S)$ for all sets of vertices $S$ and all
$\s\in \bZ_{r_1}\times\dots\times\bZ_{r_k}$, then the corresponding graph is 
democratic.

If we wish to represent the distance matrix $M^{[r]}_{ij}$ in matrix form, corresponding
to \re{pat1} or \re{M2}, then the labeling of the vertices by tuples turns out to be
a nuisance, because we would have to straighten the indices as in the $k=2$ case.
Instead, it is more natural to replace the vector space $\bR^r$ by a tensor product
$\bR^{r_1}\times\dots\times\bR^{r_k}$, with standard basis indexed by precisely
the tuples above. The generators of the symmetry group then take the form  
%We associate to this factorisation the $r_{\!\!A}{\times} r_{\!\!A}$ permutation matrices 
%p^{[A]}$, of order $r_{\!\!A}$,  
%which permute the indices $i_A$ cyclically and thus induce 
%an $r{\times} r$ permutation on the indices $i$ generated by
\be
P^{[A]}=\1^{[1]}\otimes \1^{[2]}
     \otimes\dots\otimes p^{[A]}\otimes\dots\otimes \1^{[k]}.
\end{equation}
where  $\1^{[B]}$ is the $r_{\!B}{\times} r_{\!B}$ unit matrix and
$p^{[A]}$ are the permutation matrices of order $r_{\!\!A}$,  
which permute the indices $i_A$ cyclically and thus induce 
an $r{\times} r$ permutation on the indices $i$.   
As in the $k=2$ case above, there corresponds to every factorisation 
$r=r_1r_2\cdots r_k$  a distance matrix invariant
under all the  permutations $P^{[A]}$, their powers and their products.

If $r$ has $m$ prime factors, 
$r=s_1\cdots s_m\,$, with all the $s_i$'s  different, the
number of inequivalent democratic graphs with $r$ vertices depends
only on $m$ and corresponds to the number of ways a set with $n$ 
elements can be partitioned into disjoint, non-empty subsets.
This is precisely the $m$-th Bell number $B_m$, which is given by
the formula %(see e.g. \cite{A})
\be
B_{m+1} = \sum_{k=0}^m \binom{m}{k} B_k \quad ,\quad  B_0=1\ .
\end{equation}
If some of the prime factors $s_i$ are equal, the partitions in different
subsets leading to the same set of $r_{\!A}$'s yield equivalent graphs.
We shall give further details and examples elsewhere \cite{dnw2}.

\end{document}